\documentclass[12 pt]{amsart}

\usepackage{amssymb,latexsym,amscd}
\usepackage[mathscr]{euscript}
\usepackage[all]{xy}
\usepackage{layout}

\title[Commuting elements in central products of special unitary 
groups]{Commuting elements in central products of special unitary 
groups}

\newtheorem{theorem}{Theorem}
\newtheorem*{theorem1*}{Theorem 1}
\newtheorem*{theorem2*}{Theorem 2}
\newtheorem{proposition}[theorem]{Proposition}
\newtheorem{lemma}[theorem]{Lemma}

\newtheorem{definition}[theorem]{Definition}
\newtheorem*{notation*}{Notation}

\def\Proof{\medskip\noindent{\bf Proof: }}

\def\Z{\mathbb{Z}}

\def\C{\mathbb{C}}

\def\F{\mathbb{F}}
\def\S{\mathbb{S}}

\def\M{\mathcal{M}}

\def\AC{\mathcal{AC}}
\def\uc{\underline{c}}

\def\uw{\underline{w}}
\def\ux{\underline{x}}
\def\uy{\underline{y}}
\def\uz{\underline{z}}

\def\x{\times}

\DeclareMathOperator{\Hom}{Hom}	
	
\DeclareMathOperator{\Rep}{Rep}	

\begin{document}

\author[A.~Adem]{Alejandro Adem$^{*}$}
\address{Department of Mathematics,
University of British Columbia, Vancouver BC V6T 1Z2, Canada}
\email{adem@math.ubc.ca}
\thanks{$^{*}$Partially supported by NSERC}

\author[F.~R.~Cohen]{F.~R.~Cohen$^{**}$}
\address{Department of Mathematics,
University of Rochester, Rochester NY 14627, USA}
\email{cohf@math.rochester.edu}
\thanks{$^{**}$Partially supported by 
DARPA grant number 2006-06918-01}

\author[J.~M.~G\'omez]{Jos\'e Manuel G\'omez}
\address{Department of Mathematics,
University of British Columbia, Vancouver BC V6T 1Z2, Canada}
\email{josmago@math.ubc.ca}

\begin{abstract}
In this paper the space of commuting elements in the central product
$G_{m,p}$ of $m$ copies of the special unitary group $SU(p)$ is
studied, where $p$ is a prime number. In particular, a computation
for the number of path-connected components of these spaces is given
and the geometry of the moduli space $\Rep(\mathbb Z^n, G_{m,p})$ 
of isomorphism classes of flat connections on principal 
$G_{m,p}$--bundles over the $n$--torus is completely
described for all values of $n$, $m$ and $p$.
\end{abstract}

\maketitle

\section{Introduction}

Let $G$ be a compact Lie group. The space of homomorphisms
$\Hom(\Z^{n},G)$ can be identified with the space of commuting
$n$-tuples in $G$, topologized as a subspace of the cartesian
product $G^{n}$. The quotient $\Hom(\Z^{n},G)/G$ under the conjugation 
action by $G$ is the moduli space  of isomorphism classes of flat connections 
on principal $G$-bundles over the $n$-–torus $(\S^{1})^{n}$. In the past 
few years there has been an increasing interest in understanding 
these spaces, especially in computing their number of path-connected 
components and their cohomology groups as they naturally appear in 
a number of quantum field theories such as Yang-Mills and Chern-Simons 
theories.

\medskip

In \cite{BFM} the space of commuting elements in a Lie group $G$ was
analyzed by considering the space of almost commuting elements in
the universal cover of $G$ (i.e. elements which commute up to central 
elements, see Definition \ref{Definition almost commuting}). In 
particular it was shown  that $\Rep(\Z^{n},G):=\Hom(\Z^{n},G)/G$ 
is determined by the geometry of $G$ and explicit formulations were 
given for $n=2$ and $n=3$; indeed the main focus there was to 
describe the associated moduli spaces of bundles over 
$\S^1\times\S^1$ and $\S^1\times\S^1\times\S^1$.

\medskip

On the other hand, in \cite{AC} the spaces of the form 
$\Hom(\Z^{n},G)$ were studied from a homotopical point of view. 
In particular, it was shown that if $G$ is a closed subgroup 
of $GL(n,\C)$, then there exists a natural homotopy equivalence 
after a single suspension
\begin{equation}\label{AdemCohen}
\Theta_{n}:\Sigma(\Hom(\Z^{n},G))\simeq \bigvee_{1\le r\le
n}\Sigma\left(\bigvee^{\binom{{n}}{{r}}} \Hom(\Z^{r},G)/S_{r}(G)
\right),
\end{equation}
where $S_{r}(G)\subset \Hom(\Z^{r},G)$ is the subspace of $r$-tuples
$(x_{1},...,x_{r})\in \Hom(\Z^{r},G)$ for which at least one of the
$x_{i}$ equals $1_{G}$.  In \cite{ACG}, the authors show that a
similar decomposition to (\ref{AdemCohen}) also holds for the space
of almost commuting elements in a compact Lie group $G$ and that the
corresponding map $\Theta_{n}$ is actually a $G$-equivariant
homotopy equivalence thus affording a stable decomposition for the
associated spaces of representations.

\medskip

Based on these stable homotopy equivalences it seems natural to
explore situations where the geometric description of the moduli
spaces associated to commuting pairs and triples provided in
\cite{BFM}, can be extended to arbitrary commuting $n$--tuples. In
particular it can be seen that if the maximal abelian subgroups in
$G$ are path-connected, then all of the spaces $\Hom(\Z^n,G)$ are
path-connected. However, if the fundamental group of $G$ has
$p$--torsion, then it is known (see \cite{borel2}, page 139) that
there is a subgroup $\mathbb Z/p\times\mathbb Z/p\subset G$ which is
not contained in a torus and so the spaces of commuting $n$--tuples
cannot be path--connected. Thus it is natural to consider examples
where $\pi_1(G)\cong\mathbb Z/p$.

\medskip

In this paper the spaces of the form $\Hom(\Z^{n},G_{m,p})$ are
studied, where
\[
G_{m,p}=(SU(p)^{m})/(\Delta(\Z/p))
\]
is an $m$-fold central product of $SU(p)$, for a prime $p$. Thus
these are natural examples of compact Lie groups having a
fundamental group of prime order. The study of almost commuting
elements in $SU(p)^{m}$ provides a way to compute the number of
path-connected components of $\Hom(\Z^{n},G_{m,p})$. In addition, 
the structure of the components can be explicitly described. The
following theorem summarizes these results:

\begin{theorem}\label{main theorem}
For $n\ge 1$ and $p$ a prime number, the space
$\Hom(\Z^{n},G_{m,p})$ has
\[
N(n,m,p)=\frac{p^{(m-1)(n-2)}(p^{n}-1)(p^{n-1}-1)}{p^{2}-1}+1
\]
path--connected components. The path-connected component containing 
$(1,...,1)$ is a quotient of, and has the same rational cohomology as,
$$(G_{m,p}/(\S^{1})^{m(p-1)})\x_{(\Sigma_{p})^m}(\S^{1})^{m(p-1)n},$$
whereas all the other path--connected components are homeomorphic to
$$(SU(p))^{m}/((\Z/p)^{m-1}\times E_{p}),$$
where $E_{p}\subset SU(p)$
is the quaternion group $Q_{8}$ of order eight when $p=2$ and the
extraspecial $p$--group of order $p^3$ and exponent $p$ when $p>2$.
\end{theorem}

In section \ref{commuting in G_{m,p}} it is explained how the 
path-connected component of $\Hom(\Z^{n},G_{m,p})$ containing 
$(1,...,1)$ can be seen as a quotient of the compact manifold  
\[
(G_{m,p}/(\S^{1})^{m(p-1)})\x_{(\Sigma_{p})^m}(\S^{1})^{m(p-1)n}.
\]

\medskip

A particular case of relevance of Theorem \ref{main theorem} is
the case where $m=1$. In this case, $G_{1,p}=PU(p)$ and according to
the theorem  $\Hom(\Z^{n},PU(p))$ has
\[
N(n,1,p)=\frac{(p^{n}-1)(p^{n-1}-1)}{p^{2}-1}+1
\]
path--connected components. Moreover,
$$\frac{(p^{n}-1)(p^{n-1}-1)}{p^{2}-1}$$
of these components are
homeomorphic to $SU(p)/E_{p}$. On the other hand, the number $x_n$
of path--connected components of $\Hom(\Z^{n},SO(3))$ that \textit{do not}
contain the element $(1,...,1)$ was computed in \cite{Torres}, where
it was shown that
\[
x_{n}=\left\{\begin{array}{ll}
\frac{1}{6}(4^{n}-3\times 2^{n}+2) & \mbox{if $n$ is even},\\
 & \\
\frac{2}{3}(4^{n-1}-1)-2^{n-1}+1   &\mbox{if $n$ is odd}.
\end{array}\right.
\]
Note that in Theorem \ref{main theorem} the case $p=2$ and $m=1$
corresponds to
\[
G_{1,2}=SU(2)/(\Z/2)=PU(2)\cong SO(3)
\]
which is precisely the situation already studied \cite{Torres}. It
is easy to verify that
\[
x_{n}=\frac{(2^{n}-1)(2^{n-1}-1)}{3}
\]
and thus the two approaches give the same answer. 

\medskip

Taking a quotient by the conjugation action of $G_{m,p}$ yields the
following.

\begin{theorem}\label{corollary main} The moduli space of isomorphism 
classes of flat connections on principal $G_{m,p}$--bundles over 
an $n$--torus is given by
\[
\Rep(\Z^{n},G_{m,p})\cong ((\S^{1})^{(p-1)mn}/(\Sigma_{p})^m)\sqcup
X_{n,m,p},
\]
where $X_{n,m,p}$ is a finite set with $N(n,m,p)-1$ points.
\end{theorem}

As can be expected, these quotient spaces are much simpler than the
spaces of homomorphisms lying above them, which can contain
interesting geometric information which is lost modulo conjugation;
suffice it to say that for $n=1$ this is the difference between the
group $G_{m,p}$ and its quotient under conjugation $T/W$ where
$T\subset G_{m,p}$ is a maximal torus with Weyl group $W$. Also,
it's worth noting that the components which do not correspond to the
identity element deserve special attention, as they are somewhat
exotic.

\medskip

It also seems relevant to point out that the central products
considered here arise as subgroups of some of the exceptional Lie
groups. For example
$$G_{2,2}\subset \mathbb G_2,~~~ G_{2,3}\subset \mathbb F_4,~~~ 
G_{2,5}\subset \mathbb E_8,~~~ G_{3,3}\subset \mathbb E_6$$ and they 
give rise to subgroups of the form $(\mathbb Z/p)^3$ which are not
contained in the maximal tori, thus explaining the torsion in the
cohomology of the classifying spaces of these exceptional groups
(see \cite{borel2}, pages 153--154) even though they are simply
connected. It would seem that the results here could be applied to
provide information about $\Rep(\mathbb Z^n, G)$, where $G$ is one
of these groups.

\medskip

\begin{notation*} From now on, for a prime number $p$,
$E_{p}$ denotes the quaternion group $Q_{8}$ of order eight when
$p=2$ and the extraspecial $p$--group of order $p^3$ and exponent
$p$ when $p>2$. Note that this group can be identified with the
$p$--Sylow subgroup of $SL_{3}(\F_{p})$. Also, given an integer
$m\ge 1$,
\[
G_{m,p}:=(SU(p)^{m})/(\Delta(\Z/p)),
\]
here $\Delta(\Z/p)$ is seen as a subgroup of $SU(p)^{m}$ by
considering the diagonal map 
\[
\Delta(\Z/p)\hookrightarrow (\Z/p)^{m}=Z(SU(p)^{m}).
\] 
Thus $G_{m,p}$ is the $m$-fold central
product of $SU(p)$.
\end{notation*}

\noindent{\bf Acknowledgments.} The authors would like to thank the 
referee for helpful comments and suggestions.

\section{Almost commuting elements}\label{section almost commuting}

In this section almost commuting elements in a Lie group are introduced.

\begin{definition}\label{Definition almost commuting}
Take $G$ a Lie group and $K\subset Z(G)$ a closed subgroup. An
$n$-tuple $\ux:=(x_{1},...,x_{n})\in G^{n}$ is said to be a
$K$-almost commuting $n$-tuple if $[x_{i},x_{j}]\in K\subset Z(G)$
for every $1\le i,j\le n$.
\end{definition}

The motivation for considering almost commuting elements is as
follows. Consider the space $\Hom(\Z^{n},H)$, where $H$ can be
written in the form $H=G/K$, for a Lie group $G$ and a closed
subgroup $K\subset Z(G)$. In this case, the natural map $f:G\to G/K$
is both a homomorphism and a principal $K$-bundle. If
$\ux=(x_{1},...,x_{n})$ is a sequence of elements in $G/K$ that
commute, then for any lifting $\tilde{x}_{i}$ of $x_{i}$ the
commutator $[\tilde{x}_{i},\tilde{x}_{j}]\in K\subset Z(G)$ and the
space of all such sequences can be used to study 
$\Hom(\Z^{n},G/K)$.

\begin{definition} Given a compact Lie group $G$ and 
$K\subset Z(G)$ a closed subgroup  define
\[
B_{n}(G,K)=\{(x_{1},...,x_{n})\in G^{n}~|~ [x_{i},x_{j}]
\in K \text{ for all  } i, j\}.
\]
\end{definition}

The set $B_{n}(G,K)$ can be regarded as a topological space by
naturally identifying it with a subspace of $G^{n}$.  The following
simple lemma describes the precise relationship between
$B_{n}(G,K)$ and $\Hom(\Z^{n},G/K)$.

\begin{lemma}\label{principal bundles}
Let $G$ be a Lie group and $K\subset Z(G)$ a closed subgroup. 
Then the quotient map $f:G\to G/K$
induces a $G$-equivariant principal $K^{n}$-bundle
\[
\phi_{n}:B_{n}(G,K)\to \Hom(\Z^{n},G/K).
\]
\end{lemma}

In general $K$-almost commuting elements in $G$ can be used to 
obtain a decomposition of the space $\Hom(\Z^{n},G/K)$ into 
the union of (possibly empty) open and closed subspaces in the 
following way. Given $\ux=(x_{1},...,x_{n})\in B_{n}(G,K)$ 
consider the different commutators $d_{ij}=[x_{i},x_{j}]\in K$ 
for $1\le i,j\le n$. The elements $d_{ij}$ are such that 
$d_{ii}=1$ and $d_{ij}=d_{ji}^{-1}$, thus the matrix $D=(d_{ij})$ 
is an antisymmetric matrix with entries in $K\subset Z(G)$ that 
varies continuously with $\ux$. Let $T(n,\pi_{0}(K))$ be the set 
of all $n\times n$ antisymmetric matrices $C=(c_{ij})$ with entries 
in $\pi_{0}(K)$. Given  a matrix $C\in T(n,\pi_{0}(K))$ define
\[
\AC_{G}(C)=\{(x_{1},...,x_{n})\in G^{n} ~|~
\pi_{0}([x_{i},x_{j}])=c_{ij}\in \pi_{0}(K) \}\subset B_{n}(G,K),
\]
and
\[
\Hom(\Z^{n},G/K)_{C}=\phi_{n}(\AC_{G}(C))\subset \Hom(\Z^{n},G/K).
\]

Note that both $\AC_{G}(C)$ and $\Hom(\Z^{n},G/K)_{C}$ are invariant 
under the conjugation action of $G$. Also these can be endowed with 
the natural subspace topology and in this case each 
$\Hom(\Z^{n},G/K)_{C}$ is both open and closed in 
$\Hom(\Z^{n},G/K)$ and thus a union of connected components. 
The restriction of $\phi_{n}$ defines a principal $K^{n}$-bundle
\[
\AC_{G}(C)\to \Hom(\Z^{n},G/K)_{C}
\]
and there is a decomposition
\begin{equation}\label{decomposition}
\Hom(\Z^{n},G/K)=\bigsqcup_{C\in
T(n,\pi_{0}(K))}\Hom(\Z^{n},G/K)_{C}.
\end{equation}

In \cite{BFM}, Borel, Friedman and Morgan showed that the orbit space 
$\M_{G}(C):=\AC_{G}(C)/G$ is describable in terms of the 
geometry of $G$. Moreover, they obtained explicit descriptions 
for $n=2$ and $n=3$. In the next section, their work 
will be used  to obtain an explicit description for 
$\Hom(\Z^{n},G_{m,p})$ for every $n$. This sheds some light in the 
structure of the spaces of the form $\Hom(\Z^{n},G)$ for a general 
compact Lie group $G$.

\section{Commuting elements in $G_{m,p}$}\label{commuting in G_{m,p}}

The goal of this section is to prove Theorems \ref{main theorem} 
and \ref{corollary main} in the introduction. These are the main 
results of this article and are proved using decomposition 
(\ref{decomposition}). 

\medskip

To start, suppose that $G$ is a compact connected Lie group. 
Let $\Hom(\Z^{n},G)_{(1,...,1)}$ be the path-connected 
component of $\Hom(\Z^{n},G)$ that contains $(1,...,1)$. By 
\cite[Proposition 2.3]{AC}, if every abelian subgroup of $G$ 
is contained in a path-connected abelian subgroup, then the space 
$\Hom(\Z^{n},G)$ is path-connected and thus agrees with 
$\Hom(\Z^{n},G)_{(1,...,1)}$. In \cite{Baird}, the spaces of the 
form $\Hom(\Z^{n},G)_{(1,...,1)}$ were studied. For example, 
the cohomology groups with rational coefficients of these spaces 
were computed. Some of the results proved in \cite{Baird} 
are recalled next. The reader is referred to \cite{Baird} 
for the proofs of these facts.

\medskip

Fix $T\subset G$ a maximal torus in $G$. The conjugation action 
of $G$ induces a $G$-equivariant map
\begin{align}\label{map conjugation}
\varphi_{n}:G\times T^{n}&\to \Hom(\Z^{n},G)_{(1,...,1)}\\
(g,t_{1},...,t_{n})&\mapsto (gt_{1}g^{-1},...,gt_{n}g^{-1}).
\end{align}
By \cite[Lemma 4.2]{Baird} it follows that every commuting 
$n$-tuple in $\Hom(\Z^{n},G)_{(1,...,1)}$ lies in a maximal torus 
of $G$. Since any two maximal tori in $G$ are conjugated this 
shows that the map $\varphi_{n}$ is surjective. Note that $N(T)$ 
acts on $G\times T^{n}$ diagonally and that $\varphi_{n}$ is
invariant under this action. Therefore $\varphi_{n}$ descends 
to a map
\[
\bar{\varphi}_{n}:G/T\x_{W}T^{n}=G\x_{N(T)}T^{n}\to
\Hom(\Z^{n},G)_{(1,...,1)},
\]
where $W$ is the Weyl group associated to $T$. In fact 
$G\x_{N(T)}T^{n}$ is a nonsingular real algebraic variety and 
$\bar{\varphi}_{n}$ is a resolution of singularities for 
$\Hom(\Z^{n},G)_{(1,...,1)}$ as it was pointed out in 
\cite{Baird}. Thus in general $\Hom(\Z^{n},G)_{(1,...,1)}$ is 
homeomorphic to the quotient of the compact manifold 
$G/T\times_{W} T^{n}$ where each fiber $\bar{\varphi_{n}}^{-1}(\ux)$ 
is collapsed to a point for $\ux\in \Hom(\Z^{n},G)_{(1,...,1)}$. 
Moreover, modulo the conjugation action of $G$, $\bar{\varphi_{n}}$ 
induces a homeomorphism 
\[
T^{n}/W\stackrel{\cong}{\rightarrow} \Rep(\Z^{n},G)_{(1,...,1)},
\]
with $W$ acting diagonally on $T^{n}$. In addition, by 
\cite[Theorem 4.3]{Baird} given a field $\F$ of
characteristic relatively prime to $|W|$, the map $\bar{\varphi}_{n}$
induces an isomorphism
\begin{equation}\label{Rational cohomology}
H^{*}(\Hom(\Z^{n},G)_{(1,...,1)};\F)\cong H^{*}(G/T\times
T^{n};\F)^{W}.
\end{equation}
For the case of $G=G_{m,p}$, a maximal torus $T$ is homeomorphic to
$(\S^{1})^{m(p-1)}$ and $W=(\Sigma_{p})^{m}$. Moreover, if 
$C_{1}$ is the trivial matrix whose entries are all $1$ then 
it follows that 
$\Hom(\Z^{n},G_{m,p})_{(1,...,1)}=\Hom(\Z^{n},G_{m,p})_{C_{1}}$
is a quotient of
\begin{equation}\label{aux1}
(G_{m,p}/(\S^{1})^{m(p-1)})\x_{(\Sigma_{p})^{m}}(\S^{1})^{m(p-1)n},
\end{equation}
also 
\begin{equation}\label{aux2}
\Rep(\Z^{n},G_{m,p})_{(1,...,1)}\cong 
(\S^{1})^{m(p-1)n}/(\Sigma_{p})^{m}
\end{equation}
and 
\begin{equation}\label{aux3}
H^{*}(\Hom(\Z^{n},G_{m,p})_{(1,...,1)};\F)\cong
H^*((G_{m,p}/(\S^{1})^{m(p-1)})\times
(\S^{1})^{m(p-1)n};\F)^{(\Sigma_{p})^m}.
\end{equation}
for every field $\F$ with characteristic not dividing $p!$.

\medskip

Next the spaces of the form  $\Hom(\Z^{n},G_{m,p})_{C}$ for
$C\ne C_{1}$ are studied. The following lemma, which can be 
proved directly or using \cite[Proposition 4.1.1]{BFM}, is used 
to handle this case.

\begin{lemma}\label{lemma1}
Let $c\in Z(SU(p))-\{1\}$. Then there is a pair $(x_{o},y_{o})$ of
elements in $SU(p)$ with $[x_{o},y_{o}]=c$. Moreover, the pair
$(x_{o},y_{o})$ is unique up to conjugation and if $(x,y)$ is any
such pair then $Z_{SU(p)}(x,y)=Z(SU(p))$.
\end{lemma}

The following notation will be used. Given an element 
$\uc\in \Delta(\Z/p)$, $C(\uc)$ denotes the $2\times 2$ antisymmetric 
matrix with entries in $\Delta(\Z/p)$ defined by 
$\uc_{11}=\uc_{22}=\underline{1}\in\Delta(\Z/p)$ and $\uc_{12}=
\uc_{21}^{-1}=\uc$. Theorem \ref{main theorem} will be proved by 
considering first the case $n=2$.

\begin{proposition}\label{case n=2}
The space $\Hom(\Z^{2},G_{m,p})$ has $p$ path--connected components. One
of these components is $\Hom(\Z^{2},G_{m,p})_{(1,...,1)}$ and the rest
of the components are all homeomorphic to
$SU(p)^{m}/((\Z/p)^{m-1}\times E_{p})$.
\end{proposition}
\Proof This proposition will be proved by studying the different
spaces $\AC_{SU(p)^{m}}(C)$, where $C$ is a general matrix in
$T(2,\Delta(\Z/p))$. Such a matrix is of the form $C=C(\uc)$ for
some $\uc\in \Delta(\Z/p)$. When $\uc=\underline{1}$ the space
$\AC_{SU(p)^{m}}(C(\underline{1}))$ equals $\Hom(\Z^{2},SU(p)^{m})$
which is path--connected. Thus suppose that $\uc\ne \underline{1}$. Since
$\uc\in \Delta(\Z/p)$, it is of the form $\uc=(c,...,c)$ for $c\in
\Z/p=Z(SU(p))$ with $c\ne 1$.  Fix a pair of elements $x_{o},y_{o}$
in $SU(p)$ with $[x_{o},y_{o}]=c$. By Lemma \ref{lemma1} 
the group $SU(p)^{m}$ acts transitively by conjugation on each
$\AC_{SU(p)^{m}}(C(\uc))$, thus there is a continuous surjective map
\begin{align*}
SU(p)^{m}&\to \AC_{SU(p)^{m}}(C(\uc))\\
(g_{1},...,g_{m})&\mapsto (\ux, \uy)
\end{align*}
where
\[
\ux=(g_{1}x_{o}g_{1}^{-1},...,g_{m}x_{o}g_{m}^{-1}) \text{ and }
\uy=(g_{1}y_{o}g_{1}^{-1},...,g_{m}y_{o}g_{m}^{-1}).
\]
In particular, $\AC_{SU(p)^{m}}(C(\uc))$ is path--connected and
\[
\AC_{SU(p)^{m}}(C(\uc))\cong
SU(p)^{m}/SU(p)^{m}_{(\ux_{o},\uy_{o})},
\]
where $\ux_{o}=(x_{o},...,x_{o})$, $\uy_{o}=(y_{o},...,y_{o})$ and
$SU(p)^{m}_{(\ux_{o},\uy_{o})}$ is the isotropy subgroup of
$SU(p)^{m}$ at $(\ux_{o},\uy_{o})$. Note that  
$Z_{SU(p)}(x_{o},y_{o})=Z(SU(p))$ by Lemma \ref{lemma1}, hence
\[
SU(p)^{m}_{(\ux_{o},\uy_{o})}=Z(SU(p)^{m})=\left<c\right>^{m}
=(\Z/p)^{m}
\]
and therefore
\begin{equation}\label{eqaux}
\AC_{SU(p)^{m}}(C(\uc))\cong SU(p)^{m}/\left<c\right>^{m}.
\end{equation}
On the other hand, $(\Delta(\Z/p))^{2}$ acts on $\AC_{SU(p)^{m}}(C(\uc))$ 
by left componentwise multiplication. This action  gives 
rise to a covering space sequence
\begin{equation}\label{covering PU(p)}
(\Delta(\Z/p))^{2}\to \AC_{SU(p)^{m}}(C(\uc))\to
\Hom(\Z^{2},G_{m,p})_{C(\uc)}.
\end{equation}
In particular $\Hom(\Z^{2},G_{m,p})_{C(\uc)}$ is path--connected and
\[
\Hom(\Z^{2},G_{m,p})_{C(\uc)}\cong
\AC_{SU(p)^{m}}(C(\uc))/(\Delta(\Z/p))^{2}.
\]
Notice that under the identification (\ref{eqaux}), this action of 
$(\Delta(\Z/p))^{2}$ corresponds to  
\begin{align*}
(\Delta(\Z/p))^{2}\times SU(p)^{m}/\left<c\right>^{m}
&\to SU(p)^{m}/\left<c\right>^{m}\\
(\uc^{s},\uc^{r}),[(g_{1},....,g_{m})]&\mapsto
[(g_{1}x_{o}^{r}y_{o}^{-s},...,g_{m}x_{o}^{r}y_{o}^{-s})].
\end{align*}
This is true because
\[
(x_{o}^{r}y_{o}^{-s})x_{o}(x_{o}^{r}y_{o}^{-s})^{-1}=c^{s}x_{o} 
\text{  and  }
(x_{o}^{r}y_{o}^{-s})y_{o}(x_{o}^{r}y_{o}^{-s})^{-1}=c^{r}y_{o}.
\]
It follows then that $\Hom(\Z^{2},G_{m,p})_{C(\uc)}\cong
SU(p)^{m}/K_{p}$, where $K_{p}\subset SU(p)^{m}$ is the subgroup
generated by $Z(SU(p)^{m})$, $\ux_{o}=(x_{o},...,x_{o})$ and
$\uy_{o}=(y_{o},...,y_{o})$.  By \cite[Proposition 4.1.1]{BFM} the
subgroup generated by $x_{o}$ and $y_{o}$ in $SU(p)/\left<c\right>$
is isomorphic to $(\Z/p)^{2}$ and by \cite[Corollary 4.1.2]{BFM}
$x_{o}$ and $y_{o}$ have order $4$ and $x_{o}^{2}=y_{o}^{2}=c$ if
$p=2$ and order $p$ if $p>2$. Thus, when $p=2$ the subgroup $E_{2}$
of $SU(2)$ generated by $c$, $x_{o}$ and $y_{o}$ has the
presentation
\[
E_{2}:=\left\{x,y~|~ x^{4}=y^{4}=1, x^{2}=y^{2}, yxy^{-1}=x \right\}
\]
and thus $E_{2}=Q_{8}$. When $p>2$, the subgroup $E_{p}$ of $SU(p)$
generated by $c$, $x_{o}$ and $y_{o}$ has the presentation
\[
E_{p}=\left\{x,y,c~|~ x^{p}=y^{p}=c^{p}=1, xc=cx, yc=cy, xy=cyx
\right\}
\]
and this is easily seen to be the extraspecial $p$--group
$Syl_{p}(SL_{3}(\F_{p}))$. The group $K_{p}$ fits into a short exact
sequence
\[
1\to (\Z/p)^{m-1}\to K_{p}\to \left<\uc,\ux_{o},\uy_{o} \right>\to 1,
\]
where the map $(\Z/p)^{m-1}\to K_{p}$ is as follows. Let
$u_{1},...,u_{m-1}$ be elements in the $\F_{p}$ vector space
$(\Z/p)^{m}$ such that $u_{1},...,u_{m-1},\uc$ forms a basis. Then
the $i$-th generator of $(\Z/p)^{m-1}$ is sent to $u_{i}$ for $1\le
i\le m-1$. The previous short exact sequences splits,
$\left<\uc,\ux_{o},\uy_{o} \right>\cong E_{p}$ and therefore
$K_{p}\cong (\Z/p)^{m-1}\times E_{p}$. To finish the proposition,
note there are precisely $(p-1)$ non-trivial elements 
$\uc\in \Delta(\Z/p)$.
\qed
\medskip

From the previous proposition it is deduced that $N(2,m,p)=p$. Moreover,
from the proof it follows that $G_{m,p}$ acts transitively  by 
conjugation  on each component that is homeomorphic 
to $SU(p)^{m}/((\Z/p)^{m-1}\times E_{p})$.

\begin{lemma}\label{important}
Suppose that $\ux=(x_{1},...,x_{m})$ and $\uy=(y_{1},...,y_{m})$ are
elements in $SU(p)^{m}$ that almost commute with
$\uc:=[\ux,\uy]=(c,...,c)\in \Delta(\Z/p)$ for $c\ne 1$. Take
$\uz\in SU(p)^{m}$ with $[\ux,\uz],[\uy,\uz]\in \Delta(\Z/p)$. Write
$[\ux,\uz]=\uc^{b}$ and $[\uy,\uz]=\uc^{a}$ for integers $0\le
a,b<p$. Then there is an element $\uw=(w_{1},...,w_{m})\in
Z(SU(p)^{m})$ such that $\uz=\uw \ux^{-a}\uy^{b}$; that is,
$z_{i}=w_{i}x_{i}^{-a}y_{i}^{b}$ for all $i$.
\end{lemma}
\Proof It is enough to prove the lemma for $m=1$. Fix $x,y$ and $z$ 
in $SU(p)$ such that there exists $c\in Z(SU(p))-\{1\}$ with 
$d_{1,2}:=[x,y]=c$, $d_{1,3}:=[x,z]=c^{b}$ and $d_{2,3}:=[y,z]=c^{a}$ 
for integers $0\le a,b<p$. Then the triple $(x,y,z)$ is an almost 
commuting triple in $\AC_{SU(p)}(D)$, where $D$ is the 
antisymmetric matrix with entries $d_{i,j}$. Consider the map 
\begin{align*}
\psi:\AC_{SU(p)}(C)&\to \AC_{SU(p)}(D)\\
(x_{1},x_{2},x_{3})&\mapsto (x_{1},x_{2},x_{1}^{-a}x_{2}^{b}x_{3}),
\end{align*}
where $C$ is the antisymmetric matrix with coefficients in $Z(SU(p))$
and $c_{1,2}=c_{2,1}^{-1}=c$ and $c_{i,j}=1$ else. It is 
straight--forward to check that $\psi$ is a well defined homeomorphism 
that is equivariant under the conjugation action of $SU(p)$. Let 
$(x',y',z')$ be any element in $\AC_{SU(p)}(C)$. This means that 
$[x',y']=c\ne 1$ and $z'$ commutes with both $x'$ 
and $y'$. Thus $z'\in Z_{SU(p)}(x',y')=Z(SU(p))$ 
by Lemma \ref{lemma1}. On the other hand, since $[x,y]=c$ it follows by 
Lemma \ref{lemma1} that the pair $(x',y')$ is conjugate to $(x,y)$. This 
shows that any element in $\AC_{SU(p)}(C)$ is of the form 
$(gxg^{-1},gyg^{-1},w)$ for some $g\in SU(p)$ and $w\in Z(SU(p))$. 
In particular, since $\psi$ is surjective there are
$g\in SU(p)$ and $w\in Z(SU(p))$ such that
\[
(x,y,z)=\psi(gxg^{-1},gyg^{-1},w);
\]
that is,
\[
(x,y,z)=(gxg^{-1},gyg^{-1},
g(wx^{-a}y^{b})g^{-1}).
\]
This means that $gx=xg$ and $gy=yg$, hence 
$g\in Z_{SU(p)}(x,y)=Z(SU(p))$ and therefore
\[
z=wx^{-a}y^{b}.
\]
\qed
\medskip

The next step is the proof of Theorem \ref{main theorem} in the introduction.

\begin{theorem1*}
For $n\ge 1$ and $p$ a prime number, the space
$\Hom(\Z^{n},G_{m,p})$ has
\[
N(n,m,p)=\frac{p^{(m-1)(n-2)}(p^{n}-1)(p^{n-1}-1)}{p^{2}-1}+1
\]
path--connected components. The path-connected component containing 
$(1,...,1)$ is a quotient of, and has the same rational cohomology as,
$$(G_{m,p}/(\S^{1})^{m(p-1)})\x_{(\Sigma_{p})^m}(\S^{1})^{m(p-1)n},$$
whereas all the other path--connected components are homeomorphic to
$$(SU(p))^{m}/((\Z/p)^{m-1}\times E_{p}),$$
where $E_{p}\subset SU(p)$
is the quaternion group $Q_{8}$ of order eight when $p=2$ and the
extraspecial $p$--group of order $p^3$ and exponent $p$ when $p>2$.
\end{theorem1*}

\Proof Fix $p$ a prime
number. The proof of the theorem goes by induction on $n$. For $n=1$
the theorem is trivial and for $n=2$ the theorem follows by
the Proposition \ref{case n=2}. Assume then that
$n\ge 3$. To determine the value of each $N(n,m,p)$ it will be shown
that the different $N(n,m,p)$'s satisfy the recurrence equation
\[
N(n,m,p)=p^{m-1}N(n-1,m,p)+p^{m(n-2)+n-1}-p^{m(n-2)}-p^{m-1}+1.
\]
Once this proved, by induction it follows that
\[
N(n,m,p)=\frac{p^{(m-1)(n-2)}(p^{n}-1)(p^{n-1}-1)}{p^{2}-1}+1.
\]
By (\ref{decomposition}) the space $\Hom(\Z^{n},G_{m,p})$ is a
disjoint union of the different $\Hom(\Z^{n},G_{m,p})_{C}$, where
$C$ runs through the elements in $T(n,\Delta(\Z.p))$. The different
possibilities for elements $C\in T(n,\Delta(\Z/p))$ are considered
next.
\newline

$\bullet$ {\bf{Case 1}}. Suppose that $C=C_{1}\in
T(n,\Delta(\Z/p))$ is the trivial matrix whose entries are all equal
to $1$. In this case the space $\Hom(\Z^{n},G_{m,p})_{C_{1}}$ is
path--connected and as described in (\ref{aux1}) and (\ref{aux3}).
\newline

$\bullet$ {\bf{Case 2}}. Suppose that $C\in
T(n,\Delta(\Z/p))-\{C_{1}\}$ is such that $\uc_{1,i}=1$ for all $i$.
Because $C$ is not trivial there exist $2\le i,j\le n$ such that
$\uc_{i,j}\ne 1$. Take $(\ux_{1},...,\ux_{n})\in \AC_{SU(p)^{m}}(C)$.
Since $\uc_{1,i}=1$, it follows that $\ux_{1}$ commutes with
$\ux_{i}$ for all $i$. Also, $[\ux_{i},\ux_{j}]\in
\Delta(\Z/p)-\{1\}$ and thus $\ux_{1}\in
Z_{SU(p)^{m}}(\ux_{i},\ux_{j})=Z(SU(p)^{m})$ by Lemma \ref{lemma1}.
Therefore $(\ux_{1},...,\ux_{n})\in Z(SU(p)^{m})\x
\AC_{SU(p)}(\tilde{C})$, where $\tilde{C}$ is the $(n-1)\x (n-1)$
matrix obtained from $C$ by deleting the first row and column from
$C$. In this case
\begin{align*}
\Hom(\Z^{n},G_{m,p})_{C}&= (Z(SU(p)^{m})\x 
\AC_{SU(p)^{m}}(\tilde{C}))/(\Delta(\Z/p))^{n}\\
&\cong(\Z/p)^{m}/(\Delta(\Z/p))\times \Hom(\Z^{n-1},G_{m,p})_{\tilde{C}}\\
&\cong (\Z/p)^{m-1}\times \Hom(\Z^{n-1},G_{m,p})_{\tilde{C}}.
\end{align*}
By induction each path--connected component of $\Hom(\Z^{n},G_{m,p})_{C}$
is homeomorphic to $$SU(p)^{m}/((\Z/p)^{m-1}\times E_{p})$$ with
$G_{m,p}$ acting transitively by conjugation. In addition, each
matrix $C$ of the type considered in this case determines and is
uniquely determined by the corresponding $\tilde{C}$ which is non
trivial. It follows that there are $p^{m-1}(N(n-1,m,p)-1)$ path--connected
components associated to this case.
\newline

$\bullet$ {\bf{Case 3}}. Suppose that $C\in T(n,\Delta(\Z/p))$ is
such that $\uc_{1i}\ne 1$ for some $i$. Then $2\le i\le n$ as
$\uc_{11}=\underline{1}$. Let $i$ be the smallest $i$ with
$\uc_{1i}\ne 1$, let $\uc=\uc_{1i}\in \Delta(\Z/p)$ and take
$(\ux_{1},...,\ux_{n})\in \AC_{SU(p)^{m}}(C)$. For each $2\le k\le
n$ with $k\ne i$ consider the triple $(\ux_{1},\ux_{i},\ux_{k})$.
This is an almost commuting triple with $[\ux_{1},\ux_{i}]=\uc\ne
1$. By Lemma \ref{important}, if
$\uc_{1,k}=[\ux_{1},\ux_{k}]=\uc^{b_{k}}$ and
$\uc_{i,k}=[\ux_{i},\ux_{k}]=\uc^{a_{k}}$ for integers $0\le
a_{k},b_{k}<p$, then there exist $\uw_{k}\in Z(SU(p)^{m})$ such that
$\ux_{k}=\uw_{k}\ux_{1}^{-a_{k}}\ux_{i}^{b_{k}}$ for all $k$. Note
that the integers $a_{k}$ and $b_{k}$ are uniquely determined by the
condition $0\le a_{k},b_{k}<p$ and these are in turn uniquely
determined by $\uc_{1k},\uc_{ik}$. It follows that the $n$-tuple
$(\ux_{1},...,\ux_{n})$ is uniquely determined by
$(\ux_{1},\ux_{i})$, $\uc_{1,k}, \uc_{i,k}\in \Delta(\Z/p)$ and
$\uw_{k}\in Z(SU(p)^{m})$ for $k\ne 1,i$. Moreover, if as before
$C(\uc)$ is the $2\x 2$ matrix
\begin{equation*}
C(\uc)=\left[
\begin{array}{cc}
1 & \uc \\
\uc^{-1} & 1
\end{array}%
\right],
\end{equation*}
then the map
\begin{align*}
\psi:\AC_{SU(p)^{m}}(C(\uc))\x (Z(SU(p)^{m}))^{n-2}&\to 
\AC_{SU(p)^{m}}(C)\\
((\ux_{1},\ux_{i}),(\uw_{2},...,\uw_{i-1},\uw_{i+1},...,\uw_{n}))&
\mapsto(\uy_{1},...,\uy_{n}),
\end{align*}
is a homeomorphism where
\begin{equation*}
\uy_{k}=\left\{
\begin{array}{cl}
\ux_{1} & \text{if }k=1, \\
\ux_{i} & \text{if }k=i, \\
\uw_{k}\ux_{1}^{-a_{k}}\ux_{i}^{b_{k}} & \text{if }k\neq 1,i.%
\end{array}%
\right.
\end{equation*}
The map $\psi$ is $SU(p)^{m}$-equivariant, with $SU(p)^{m}$ acting
by conjugation. By passing to the quotient of the respective
$(\Delta(\Z/p))^{n}$--actions, it follows that $\psi$ induces a
homeomorphism
\[
\Hom(\Z^{2},G_{m,p})_{C(\uc)}\times (\Z/p)^{(m-1)(n-2)}\to
\Hom(\Z^{n},G_{m,p})_{C}.
\]
By the case $n=2$, each path--connected component of
$\Hom(\Z^{n},G_{m,p})_{C}$ is of the desired type and there are
$p^{(m-1)(n-2)}$ such components associated to $C$. It also follows
that $G_{m,p}$ acts transitively on each of these components.
Moreover, $C$ is uniquely determined by $\uc=\uc_{1i}\ne 1$,
$\uc_{1k}$ for $i+1\le k\le n$ and $\uc_{ik}$ for $2\le k\le n$ and
$k\ne i$. Thus there are in total $p^{(m-1)(n-2)}(p-1)p^{2n-i-2}$
different components associated to such $C$ with $\uc_{1,i}\ne 1$.
Letting $2\le i\le n$ vary, a total number of
\[
\sum_{i=2}^{n}p^{(m-1)(n-2)}(p-1)p^{2n-i-2}=p^{m(n-2)}(p^{n-1}-1)
\]
path--connected components is obtained for this case.

\medskip

By adding the contributions from case 1, case 2 and case 3 the
recurrence equation
\begin{align*}
N(n,m,p)&= 1+p^{m-1}(N(n-1,m,p)-1)+ p^{m(n-2)}(p^{n-1}-1)\\
&=p^{m-1}N(n-1,m,p)+p^{m(n-2)+n-1}-p^{m(n-2)}-p^{m-1}+1.
\end{align*}
is obtained as claimed. \qed

\medskip

As mentioned before $G_{m,p}$ acts transitively on the components of 
$\Hom(\Z^{n},G_{m,p})$ that are homeomorphic to 
$SU(p)^{m}/(\Z/p^{m-1}\times E_{p})$. This shows that these 
path-connected components represent isolated points in the moduli 
space $\Rep(\Z^{n},G_{m,p})$. On the other hand, 
by (\ref{aux2}) there is a homeomorphism 
\[
\Rep(\Z^{n},G_{m,p})_{(1,...,1)}\cong 
(\S^{1})^{m(p-1)n}/(\Sigma_{p})^{m}.
\]
As a corollary of this the following theorem is obtained.

\begin{theorem2*}
Let $p$ be a prime number and $m\ge 1$. Then $\Rep(\Z^{n},G_{m,p})$
has
\[
N(n,m,p)=\frac{p^{(m-1)(n-2)}(p^{n}-1)(p^{n-1}-1)}{p^{2}-1}+1
\]
path--connected components and
\[
\Rep(\Z^{n},G_{m,p})\cong ((\S^{1})^{(p-1)mn}/(\Sigma_{p})^m)\sqcup
X_{n,m,p}.
\]
where $X_{n,m,p}$ is a finite set with $N(n,m,p)-1$ points.
\end{theorem2*}

The component of the identity can be described more explicitly as
follows. $\Sigma_p$ acts on $(\S^1)^p$ as the Weyl group of a
maximal torus in $SU(p)$. Then the product $(\Sigma_p)^m$ acts on
the product $(\S^1)^{(p-1)m}$, and therefore diagonally on the
product $(\S^1)^{(p-1)mn}$. For example, if $p=2$, the action of
$(\Sigma_2)^m$ on $(\S^1)^m$ is simply given as a product of the
complex conjugation action, and this is extended to a diagonal
action on $((\S^1)^m)^n$.

\medskip

In Theorem \ref{main theorem} if $p$ is not longer assumed to be 
a prime number then the situation is more complicated. For 
example when $n=2$, the conjugation action of $SU(r)$ on 
$\AC_{SU(r)}(C(c))$ is not longer transitive, unless $c$
is a generator of the the cyclic group $\Z/r$. Because of this, 
the space $\Hom(\Z^{n},SU(r))$ has in general more path--connected 
components that have orbifold singularities. In particular, for 
$n=2$ there is the following proposition that can be proved in 
the same way as Lemma \ref{case n=2}.

\begin{proposition}
The space $\Hom(\Z^{2},PU(r))$ has $r$ path--connected components. 
Of these $\varphi(r)$ are homeomorphic to $PU(r)/(\Z/r)^{2}$.
\end{proposition}


\begin{thebibliography}{99}
\bibitem{AC} A. Adem and F. R. Cohen. Commuting elements 
and spaces of homomorphisms.  Math. Ann.  338  (2007),  
no. 3, 587-626.

\bibitem{ACG} A. Adem, F. R. Cohen and J.M G\'omez. Stable splittings, 
spaces of representations and almost commuting elements in Lie groups. 
To appear in Math. Proc. Camb. Phil. Soc.

\bibitem{Baird} T. J. Baird. Cohomology of the space of commuting 
$n$-tuples in a compact Lie group. Algebr. Geom. Topol. 7 (2007), 
737--754.

\bibitem{borel2} A. Borel. \textbf{Collected Papers, Vol. II}. 
Springer--Verlag (1983).

\bibitem{BFM} A. Borel, R. Friedman and J. W. Morgan. Almost commuting 
elements in compact Lie groups. Mem. Amer. Math. Soc. 157 (2002), no. 747.

\bibitem{Torres} E. Torres-Giese and D. Sjerve, Fundamental groups 
of commuting elements in Lie groups. Bulleting of the London Mathematical 
Society 2008 40(1) 65-76.


\end{thebibliography}
\end{document}